\documentclass[]{amsart}
\usepackage[utf8]{inputenc}
\usepackage[mathscr]{eucal} 
\usepackage{stmaryrd} 
\usepackage{amsmath,amsthm,amsfonts,amssymb,amscd} 
\usepackage[dvipsnames]{xcolor} 
\usepackage{xspace} 
\usepackage{fancyhdr} 
\usepackage{graphicx} 
\usepackage{listings} 
\usepackage[]{hyperref} 
\usepackage{enumitem} 
\usepackage{relsize} 
\usepackage{tikz-cd} 
\usepackage{mdwlist} 
\usepackage{multicol} 
\usepackage{float} 
\usepackage{adjustbox} 
\usepackage{tikz} 
\usepackage[noabbrev,nameinlink]{cleveref} 
\Crefformat{section}{#2\S#1#3}
\usepackage{setspace} 
\usepackage{bbm} 
\usetikzlibrary{matrix, calc, arrows}

\hypersetup{
    colorlinks=true, 
    linktoc=all,     
    linkcolor=Brown,citecolor=Brown,urlcolor=MidnightBlue,  
}

\setlist[enumerate]{label=(\alph*)}

\usetikzlibrary{graphs,decorations.pathmorphing,decorations.markings}
\tikzcdset{scale cd/.style={every label/.append style={scale=#1},
        cells={nodes={scale=#1}}}}

\newcommand{\Z}{\mathbb{Z}}

\newcommand{\Q}{\mathbb{Q}}
\newcommand{\N}{\mathbb{N}}
\newcommand{\F}{\mathbb{F}}

\newcommand{\Hom}{\operatorname{Hom}}
\newcommand{\Aut}{\operatorname{Aut}}

\newcommand{\tr}{\operatorname{tr}}

\newcommand{\inv}{^{-1}}

\newcommand{\End}{\operatorname{End}}

\newcommand{\Res}{\operatorname{Res}}
\newcommand{\Ind}{\operatorname{Ind}}

\newcommand{\Syl}{\operatorname{Syl}}
\newcommand{\GL}{\operatorname{GL}}
\newcommand{\triv}{\mathbf{triv}}

\newcommand{\calA}{\mathcal{A}}

\newcommand{\calE}{\mathcal{E}}
\newcommand{\calF}{\mathcal{F}}

\newcommand{\calO}{\mathcal{O}}

\newcommand{\calT}{\mathcal{T}}

\newcommand{\calX}{\mathcal{X}}

\newcommand{\Gal}{\operatorname{Gal}}

\newcommand{\SL}{\operatorname{SL}}
\newcommand{\catmod}{\mathbf{mod}}

\newcommand{\Br}{\operatorname{Br}}

\newtheorem{theorem}{Theorem}[section]
\newtheorem{lemma}[theorem]{Lemma}
\newtheorem{prop}[theorem]{Proposition}
\newtheorem{corollary}[theorem]{Corollary}

\newtheorem*{theorem*}{Theorem}
\newtheorem{conjecture}[theorem]{Conjecture}

\theoremstyle{remark}
\newtheorem{remark}[theorem]{Remark}

\newtheorem{convention}[theorem]{Convention}
\newtheorem{observation}[theorem]{Observation}

\theoremstyle{definition}
\newtheorem{definition}[theorem]{Definition}
\newtheorem{notation}[theorem]{Notation}
\newtheorem{construction}[theorem]{Construction}

\AddToHook{env/corollary/begin}{\crefalias{theorem}{corollary}}
\AddToHook{env/proposition/begin}{\crefalias{theorem}{proposition}}
\AddToHook{env/lemma/begin}{\crefalias{theorem}{lemma}}
\AddToHook{env/definition/begin}{\crefalias{theorem}{definition}}
\AddToHook{env/example/begin}{\crefalias{theorem}{example}}
\AddToHook{env/remark/begin}{\crefalias{theorem}{remark}}
\AddToHook{env/observation/begin}{\crefalias{theorem}{observation}}

\begin{document}
    \title{On endosplit $p$-permutation resolutions and Brou\'{e}'s conjecture for $p$-solvable groups}
    \author{Sam K. Miller}
    \address{Department of Mathematics, University of Georgia, Athens, GA 30601} 
    \email{sam.miller@uga.edu} 
    \begin{abstract}
        Endosplit $p$-permutation resolutions play an instrumental role in verifying Brou\'{e}'s abelian defect group conjecture in numerous cases. We give a new characterization of all endosplit $p$-permutation resolutions and reduce the question of Galois descent of an endosplit $p$-permutation resolution to the Galois descent of the module it resolves. This is shown using techniques from the study of endotrivial complexes, the invertible objects of the bounded homotopy category of $p$-permutation modules. As an application, we show that a refinement of Brou\'{e}'s conjecture proposed by Kessar--Linckelmann holds for certain blocks of groups $G$ satisfying $G = O_{p',p,p'}(G)$ with abelian Sylow $p$-subgroup, the key reduction step in Harris--Linckelmann's verification of Brou\'e's conjecture for all $p$-solvable groups.
    \end{abstract}

    \maketitle
    \markleft{\MakeUppercase{On endosplit $p$-permutation resolutions and Brou\'{e}'s conjecture}}
    \markright{\MakeUppercase{On endosplit $p$-permutation resolutions and Brou\'{e}'s conjecture}}
    \section{Introduction}

    If one can make a conjecture about blocks, Galois-invariance should be involved too. This philosophy motivates numerous recent refinements of the Alperin--McKay conjecture proposed by Isaacs--Navarro in \cite{IN02}, Navarro in \cite{N04}, and Turull in \cite{T13}, a refinement of Alperin's weight conjecture proposed by Navarro in \cite{N04}, and a recent refinement of Brou\'{e}'s abelian defect group conjecture proposed by Kessar--Linckelmann in \cite{KL18}. The refinements of the Alperin--McKay conjecture add additional structure to the proposed bijections by imposing additional Galois-invariance or structural conditions. In particular, these refinements suggest that the conjecture should hold over non-splitting fields, descending to equivalences over finite fields of characteristic $p$ of any size. Kessar--Linckelmann's refinement of the abelian defect group conjecture predicts exactly this: it predicts that the abelian defect group conjecture should hold over \textit{any} choice of $p$-modular system, in particular, the $p$-modular system $(\Q_p, \Z_p, \F_p).$ Their refinement has been verified in numerous cases, including blocks with cyclic defect groups in \cite{KL18}, blocks with Klein-four defect groups in \cite{H23}, blocks of alternating groups in \cite{H24c}, blocks of $\SL_2(q)$ and $\GL_2(q)$ in \cite{HLZ23}, and unipotent blocks of $\GL_n(q)$ in \cite{HLZ23b}. The refined conjecture also holds for blocks of symmetric groups, as the original proof provided in \cite{CR08} holds independent of choice of field.

    The goal of this paper is twofold. Our first minor objective is to provide additional evidence that Kessar--Linckelmann's refinement of Brou\'{e}'s conjecture is true by showing it holds for particular blocks of $p$-solvable groups $G$, satisfying $G = O_{p',p,p'}(G)$, with abelian Sylow $p$-subgroup. The case we consider is the key reduction step in Harris--Linckelmann's verification of the conjecture for all blocks of $p$-solvable groups \cite{HL00}. This is a minor improvement in the literature, a step up from the result of \cite{M24d} in which we verified the conjecture for blocks of $p$-nilpotent groups with abelian Sylow $p$-subgroups. However, if one can prove an analogue of \cite[Proposition 3.1]{HL00} over arbitrary fields of characteristic $p$, a verification of Kessar--Linckelmann's conjecture for all $p$-solvable groups could follow.

    \begin{theorem}{(Theorem \ref{thm:mainthm})}
        Let $G$ be a $p$-solvable group satisfying $G = O_{p',p,p'}(G)$ with a abelian Sylow $p$-subgroup $P$. Let $b'$ be a $G$-stable block of $\overline{\F}_pO_{p'}(G)$, and let $b$ be a block idempotent of $\F_p G$ for which $bb' \neq 0$. Let $c$ be the Brauer correspondent of $b$ in $\F_pN_G(P)$. Then $b$ and $c$ are splendidly Rickard equivalent.
    \end{theorem}

    Harris--Linckelmann's construction relies on the well-understood block theory of $p$-solvable groups, and uses the machinery of \textit{endosplit $p$-permutation resolutions}.

    \begin{definition}{\cite{R96}}
        Let $G$ be a finite group. An \textit{endosplit $p$-permutation resolution} of a finitely generated $kG$-module $M$ is a bounded complex $C$ of $p$-permutation $kG$-modules together with an isomorphism between $M$ and the degree zero homology of $C$ such that
        \begin{enumerate}
            \item the homology of $C$ is concentrated in degree zero, and
            \item the complex $\End_k(C)$ is split as a complex of $kG$-modules.
        \end{enumerate}
    \end{definition}

    Such resolutions were introduced by Rickard in his verification of the abelian defect group conjecture for $p$-nilpotent groups, and have since become objects of independent interest. Our second objective, and main achievement of this paper, is to give a complete characterization of all endosplit $p$-permutation resolutions. Moreover, we prove that the question of Galois descent of endosplit $p$-permutation resolutions reduces to the question of Galois descent of the modules they resolve. If $k'/k$ is a field extension, we say a $k'G$-module $M'$ \textit{descends to $k$} if there exists a $kG$-module $M$ such that $k' \otimes_k M \cong M'$, and similarly for chain complexes.

    \begin{theorem}{(Corollary \ref{cor:dadechar} and Theorem \ref{thm:descentofepprs})}
        \begin{enumerate}
            \item Let $C$ be an endosplit $p$-permutation resolution with contractible no direct summands of a $kG$-module. Then $C$ is a direct summand of the endosplit $p$-permutation resolution $\bigoplus_{i=1}^l \Ind^G_{P_i} \Res^G_{P_i} Y$ for some \emph{chain tensor of relative syzygies} $Y$ (see Definition \ref{not:tensorsofrelsyz}) and family of $p$-subgroups $\{P_i\}_{i=1}^l$ of $G$.
            \item Let $k'/k$ be an extension of finite fields. Let $C$ be an endosplit $p$-permutation resolution for an endo-$p$-permutation $k'G$-module $M$ which contains no contractible summands. Then $C$ descends to $k$ if and only if $M$ does.
        \end{enumerate}
    \end{theorem}

    We note that normally, relative syzygies refer to modules that arise as the nonzero homology of chain complexes. In this paper, by a relative syzygy, we mean the corresponding chain complex.
    Our description is heavily inspired by Dade's description of all endopermutation $kP$-modules for finite $p$-groups $P$ in terms of the \emph{capped} endopermutation $kP$-modules \cite{D78a}. In the module-theoretic case, the question of classifying the capped endopermutation modules is a thoroughly challenging task, and was closed out nearly 30 years later by Bouc \cite{Bo06}. However in our case, tensor products of relative syzygies are easily described as products of certain two-term chain complexes.

    The paper is organized as follows. Section 2 recalls some prerequisite information on projectivity relative to a module and relatively endotrivial modules. Section 3 concerns endosplit $p$-permutation resolutions and relatively endotrivial complexes, and builds up to the classification of endosplit $p$-permutation resolutions. Section 4 concerns the Galois descent of endosplit $p$-permutation resolutions. Finally, Section 5 applies the results of the previous sections towards verifying Kessar--Linckelmann's refinement of Brou\'{e}'s conjecture.

    \textbf{Notation and conventions:} Throughout this paper, $G$ will be a finite group, $p$ a prime, and $(K, \calO, k)$ will denote a $p$-modular system. That is, $\calO$ is a complete discrete valuation ring, $k$ is its residue field, and $K$ is its field of fractions. Denote by $s_p(G)$ the set of $p$-subgroups of $G$, and $\Syl_p(G)$ the set of Sylow $p$-subgroups of $G$. The additive group of superclass functions valued on $p$-subgroups $s_p(G) \to \Z$ is denoted $CF(G,p)$.

    All modules and algebras are assumed finitely generated, and all chain complexes are assumed bounded. We denote the category of finite $G$-sets by ${}_G\mathbf{set}$, the category of finitely generated $kG$-modules by ${}_{kG}\catmod$, and the category of finitely generated $p$-permutation $kG$-modules by ${}_{kG}\triv$. Given an additive category $\calA$, $Ch^b(\calA)$ denotes the bounded chain complex category of $\calA$, and $K^b(\calA)$ denotes the bounded homotopy category of $\calA$. The symbol $\simeq$ denotes homotopy equivalence of chain complexes, i.e., isomorphism in $K^b(\calA)$, and given a $kG$-module $M$, we write $M[i]$ to denote the chain complex with $M$ in degree $i$ and zero in all other degrees.

    Given any $P$-subgroup $P$ of $G$, the Brauer construction is denoted as a functor $-(P)\colon {}_{kG}\triv\to {}_{k[N_G(P)/P]}\triv$ or $Ch^b({}_{kG}\triv) \to Ch^b({}_{k[N_G(P)/P]}\triv)$. The Brauer homomorphism is denoted $\Br_P\colon (kG)^P \to kC_G(P)$.

    For this paper, we assume familiarity with splendid Rickard equivalences, $p$-permutation modules, the Brauer construction, and some standard block theory and homological algebra; we refer the reader to \cite{L181} and \cite{L182} for a detailed treatment of these topics. We refer the reader to \cite[Section 12]{Bou10} for an overview of endopermutation modules for $p$-groups and the Dade group.

    \section{Relative $V$-projectivity and relatively endotrivial modules}

    For this section, $V$ is a $kG$-module, possibly 0. We review the notion of projectivity relative to a module, which was first considered by Okuyama in an unpublished manuscript \cite{O}, and relative endotriviality, which was first introduced by Lassueur in \cite{CL11}. We also introduce the generalized Dade group, introduced by Lassueur in \cite{CL13}, which will be relevant in the sequel.

    \begin{definition}
        Let $M$ be a $kG$-module.
        \begin{enumerate}
            \item $M$ is \textit{relatively $V$-projective} if there exists a $kG$-module $N$ for which $M \mid V\otimes_k N$. Denote the full subcategory of ${}_{kG}\catmod$ consisting of $V$-projective modules by ${}_{kG}\catmod(V),$ and denote by ${}_{kG}\triv(V)$ the full subcategory consisting of $p$-permutation $V$-projective modules.
            \item \cite[Definition 3.1.1]{CL11} $M$ is \textit{relatively $V$-endotrivial} if $M^* \otimes_k M \cong k \oplus N$ for some relatively $V$-projective $kG$-module $N$.
            \item \cite[Definition 5.3]{CL13} Set $V(\calF_G) := \bigoplus_{P \in s_p(G)\setminus \Syl_p(G)} k[G/P]$. $M$ is a \textit{strongly capped endo-$p$-permutation $kG$-module} if $M$ is endo-$p$-permutation, i.e. $M^* \otimes_k M$ is $p$-permutation, and $M$ is $V(\calF_G)$-endotrivial. If $G$ is a $p$-group, one may easily verify $M$ is strongly capped endo-$p$-permutation if and only if $M$ is capped endopermutation, i.e. $M^*\otimes_k M$ is a permutation module with the trivial module as a direct summand.
        \end{enumerate}

        If $V = 0$, we adapt the convention that the zero module is the only $V$-projective.
    \end{definition}

    \begin{remark}
        The notion of projectivity relative to a module generalizes projectivity relative to a family of subgroups. By Frobenius reciprocity, one may show that a $kG$-module $M$ is $H$-projective for some $H \leq G$ if and only if $M$ is $k[G/H]$-projective. Therefore, if $\calX$ is a subset of the set of subgroups of $G$, then a $kG$-module $M$ is projective relative to $\calX$ if and only if it is projective relative to $\bigoplus_{H \in \calX}k[G/H]$.

        Conversely, if $V$ is $p$-permutation, let $\calX_V$ denote the set of all $P \in s_p(G)$ for which $V(P) \neq 0$. Then projectivity relative to $V$ is equivalent to projectivity relative to $\calX_V$; see \cite[Theorem 3.7(c)]{M24b}. In particular, if $\calX \subseteq s_p(G)\setminus \Syl_p(G)$, then $\calX$-projectivity implies $V(\calF_G)$-projectivity.
    \end{remark}

    Projectivity relative to a module is only interesting in some cases, which motivates the next definition.

    \begin{definition}
        The $kG$-module $V$ is \textit{absolutely $p$-divisible} if every direct summand of $V$ has $k$-dimension divisible by $p$. We adapt the convention that the zero module is absolutely $p$-divisible.
    \end{definition}

    \begin{prop}{\cite[Propositions 2.2.2 and 3.4.1(a)]{CL11}}
        Let $k$ be an algebraically closed field. The following are equivalent.
        \begin{enumerate}
            \item The trivial $kG$-module $k$ is not $V$-projective.
            \item $V$ is absolutely $p$-divisible.
            \item ${}_{kG}\catmod(V) \subsetneq {}_{kG}\catmod$.
        \end{enumerate}
    \end{prop}

    \begin{convention}
        For this reason, relative $V$-projectivity is interesting only when $V$ is absolutely $p$-divisible. We remark that in the previous proposition, $k$ is assumed to be algebraically closed in \cite{CL11}. However, this will not pose a problem in the sequel, as we will primarily on relative projectivity with respect to the module $V(\calF_G)$, and it is straightforward to verify that all three of the above conditions hold for this module, regardless of algebraic closure.
    \end{convention}

    Similarly to the regular endotrivial module setting, every relatively endotrivial module has a unique indecomposable endotrivial summand.

    \begin{prop}{\cite[Lemma 3.3.1]{CL11}}
        Let $M$ be a relatively $V$-endotrivial $kG$-module with $V$ absolutely $p$-divisible. $M \cong M_0 \oplus U$, where $M_0$ is an indecomposable relatively $V$-endotrivial $kG$-module and $U$ is a $V$-projective $kG$-module or 0. Moreover, $M_0$ has vertex set $\Syl_p(G)$.
    \end{prop}

    \begin{construction}
        Let $V$ be an absolutely $p$-divisible $kG$-module.
        \begin{enumerate}
            \item \cite[Proposition 3.5.1]{CL11} We impose an equivalence relation $\sim_V$ on the class of $V$-endotrivial $kG$-modules as follows: if $M$ and $N$ are $V$-endotrivial, write $M \sim_V N$ if and only if $M$ and $N$ have isomorphic indecomposable $V$-endotrivial summands. Let $T_{k,V}(G)$ denote the resulting set of equivalence classes. $T_{k,V}(G)$ is an abelian group with addition defined by $[M] + [N] := [M\otimes_k N]$. We have $-[M] = [M^*]$ and each equivalence class has a unique indecomposable representative.

            \item \cite[Definition 5.5]{CL13} Let $D_k(G)$ denote the set of equivalence classes of strongly capped endo-$p$-permutation $kG$-modules under the equivalence relation $\sim_{V(\calF_G)}$. As before, $D_k(G)$ forms an abelian group with addition defined by $[M] + [N] := [M\otimes_k N]$. Moreover, $D_k(G)$ identifies as a subgroup of $T_{k,V(\calF_G)}(G)$. $D_k(G)$ is the \textit{generalized Dade group} of $G$. If $G$ is a $p$-group, $D_k(G)$ we recover the classical Dade group $D_k(G)$. We refer the reader to \cite[Section 12]{Bou10} for a comprehensive treatment of the classical Dade group.

            \item Given any $G$-set $X$, write $\Delta(X)$ for the $kG$-module given by the kernel of the augmentation homomorphism $kX \to k$. Let $S \in \Syl_p(G)$. If $X^S\neq \emptyset $, then $\Delta(X)$ is an endo-$p$-permutation module, see \cite{Al01}. We write $\Omega_X$ for the equivalent class of $\Delta(X)$ in $D_k(G)$.

            Extending notation, write $\Delta(X)\inv$ for the $kG$-module given by the cokernel of the inclusion $k \to kX$. We have $[\Delta(X)\inv] = -\Omega_X\in D_k(G)$. We define $D^\Omega(G)$ to be the subgroup of $D_k(G)$ generated by elements of the form $\Omega_X$, where $X$ runs over all $G$-sets $X$ satisfying $X^S = \emptyset$. In fact, $D^\Omega(G)$ is generated by elements of the form $\Omega_{G/P}$, where $P$ runs over a set of conjugacy class of representatives of $s_p(G) \setminus \Syl_p(G)$, see \cite[Lemma 12.1]{CL13}. Note that $D^\Omega(G)$ does not depend on the size of the base field $k$. For $P \in s_p(G)\setminus \Syl_p(G)$, we write $\Omega_P := \Omega_{G/P}$ for shorthand.
        \end{enumerate}

    \end{construction}

    The next proposition asserts that for $p$-permutation modules, the only types of relative projectivity in fact arise from projectivity relative to subgroups. Given a subset $\calX$ of the set of subgroups of $G$, we write ${}_{kG}\catmod(\calX)$ to denote the full subcategory of $\calX$-projective modules, and similarly for ${}_{kG}\triv(\calX)$.

    \begin{prop}{\cite[Theorem 3.7(c)]{M24b}}
        Let $V$ be a $kG$-module. Then, there exists a $p$-permutation module $W$ such that ${}_{kG} \triv(V) = {}_{kG}\triv(W)$. In particular, there exists a poset of $p$-subgroups $\calX\subseteq s_p(G)$ closed under $G$-conjugation and taking subgroups for which ${}_{kG}\triv(V) = {}_{kG}\triv(\calX)$.
    \end{prop}
    \begin{proof}
        First we show that if $M$ is a trivial source module with vertex $P \in s_p(G)$ which belongs to ${}_{kG}\triv(V)$, then $k[G/P]$ belongs to ${}_{kG}\triv(V)$. If $M \mid V\otimes_k N$ for some $kG$-module $N$, then $\Res^G_S M \mid \Res^G_S (V\otimes_k N)$ for $S \in \Syl_p(G)$. Since $M$ has vertex $P$, $\Res^G_S M$ has a direct summand isomorphic to $k[S/P]$, and any other nonisomorphic direct summand is isomorphic to $k[S/Q]$ for some subgroup $Q$ $G$-conjugate to a subgroup of $P$. Therefore $\Res^G_S (V\otimes_k N)$ has a direct summand isomorphic to $k[S/P]$. We have $\Res^G_S (V\otimes_k N) \cong  k[S/P] \oplus N'$, for some $kS$-module $N$, therefore, $k[S/P]$ is $\Res^G_S V$-projective. Therefore by \cite[Proposition 2.2.2]{CL12}, ${}_{kS}\triv(P) \subseteq {}_{kS}\triv(\Res^G_S V)$. This holds for all $G$-conjugates of $P$.

        We have that \[\Res^G_S k[G/P] \cong \bigoplus_{x \in [S\backslash G/P]} k[S/S\cap {}^xP],\] and for any $x \in G$, $k[S/S\cap {}^xP]$ is ${}^x P$-projective. Therefore, \[\Res^G_S k[G/P] \in {}_{kS}\triv(P)\cup \cdots \cup{}_{kS}\triv({}^xP) ,\] where the sum ranges over all $G$-conjugacy classes of $P$. Hence $\Res^G_S k[G/P] \in {}_{kS}\triv(\Res^G_S V)$ as well. Equivalently, $k[G/P]$ belongs to ${}_{kG}\catmod(V)$, as desired.

        Now, let $\calX$ be the set of all vertices which occur for all trivial source modules belonging to ${}_{kG}\triv(V)$. The previous claim and \cite[Proposition 2.2.2]{CL12} imply that \[{}_{kG}\triv(V) \supset {}_{kG}\triv\left(\bigoplus_{P \in \calX} k[G/P]\right) = {}_{kG}\triv(\calX).\]
        The converse inclusion follows from \cite[Proposition 3.5]{M24b}.
    \end{proof}

    Therefore, when we want to consider ${}_{kG}\triv(V)$, it suffices to assume $V$ is $p$-permutation, or equivalently consider ${}_{kG}\triv(\calX_V)$, where $\calX_V$ is the set of $p$-subgroups $P$ of $G$ for which $V(P) \neq 0$.

    \section{Endosplit $p$-permutation resolutions and relatively endotrivial complexes}

    For this section, $V$ is an \textit{absolutely $p$-divisible $kG$-module}, possibly 0. We review the notions of endosplit $p$-permutation resolutions, first introduced by Rickard in \cite{R96}, a class of ``invertible'' endotrivial $p$-permutation resolutions, endotrivial complexes, introduced in \cite{M24a}, and their relative counterpart, introduced in \cite{M24b}.

    \begin{definition}{\cite[Section 7]{R96}}
        Let $M$ be a $kG$-module. An \textit{endosplit $p$-permutation resolution} of $M$ is a bounded chain complex $C_M$ of $p$-permutation $kG$-modules such that the homology of $C_M$ is concentrated in a single degree $i \in \Z$, $H_i(C_M) \cong M$, and the complex $\End_k(C_M) \cong C_M^* \otimes_k C_M$ is split as a chain complex of $kG$-modules.
    \end{definition}

    \begin{remark}
        It is straightforward to see via the K\"unneth formula that a bounded chain complex $C$ of $p$-permutation modules is an endosplit $p$-permutation resolution (for some $kG$-module) if and only if $C^* \otimes_k C \simeq N[0]$ for some $kG$-module $N$, since split complexes are homotopy equivalent to their homology regarded as a chain complex with zero differential, and the K\"unneth formula implies $N \cong H_i(C) \otimes_k H_i(C)^*$ for a unique integer $i \in \Z$. Observe that since $N[0]$ a direct summand of $C^* \otimes_k C$, a complex of $p$-permutation modules, $N$ is $p$-permutation, so $H_i(C)$ is endo-$p$-permutation. In other words, if a $kG$-module has an endosplit $p$-permutation resolution, it is endo-$p$-permutation.

        In the original definition of an endosplit $p$-permutation resolutions, Rickard specifies that the homology of $C_M$ should be concentrated in degree zero. We adopt a looser definition, allowing the nonzero homology to be concentrated in any degree for consistency later on.
    \end{remark}

    \begin{prop}\label{prop:equivendosplit}
        Let $C$ be a bounded chain complex of $p$-permutation $kG$-modules. The following are equivalent:
        \begin{enumerate}
            \item $C$ is an endosplit $p$-permutation resolution.
            \item For all $P \in s_p(G)$, one of the following holds. Either $C(P)$ is contractible or the nonzero homology of $C(P)$ is concentrated in exactly one degree.
        \end{enumerate}
    \end{prop}
    \begin{proof}
        First, suppose $C$ is an endosplit $p$-permutation resolution for some $kG$-module $M$. Fix $P \in s_p(G)$. $\End_k(C) \cong C^* \otimes_k C$ is split, and the K\"unneth formula implies $C^* \otimes_k C$ has nonzero homology concentrated in degree zero, with $H_0(C^* \otimes_k C) \cong M^* \otimes_k M$, therefore \[C^* \otimes_k C\simeq (M^* \otimes_k M)[0].\] Applying the Brauer quotient, we obtain \[C(P)^* \otimes_k C(P) \cong (C^* \otimes_k C)(P) \simeq (M^* \otimes_k M)(P)[0].\] If $(M^* \otimes_k M)(P) \neq 0$, it follows again by the K\"unneth formula that $C(P)$ has nonzero homology in exactly one degree. Otherwise, if $(M^* \otimes_k M)(P) = 0$, then $C(P)^*\otimes_k C(P)\cong \End_k(C(P))$ is contractible. Thus $\End_{K^b({}_{k[N_G(P)/P]}\triv)}(C(P)) = 0$, and $C(P)$ is zero in $K^b({}_{kG}\triv)$, as desired.

        Now, suppose condition (b) holds. Suppose for contradiction that $C^* \otimes_k C \not\simeq M[0]$, then we may write $C^* \otimes_k C \simeq D$ where $D$ is assumed to have no contractible summands. Since $C$ is bounded, there exists a unique positive natural number $i > 0$ for which $D_i \neq 0$ and $D_j = 0$ for all $j > i$. Since $C(P)$ is either contractible or has homology contained in exactly one degree, $C(P)^* \otimes_k C(P)$ has nonzero homology concentrated in degree 0, possibly also zero. The same holds for $D(P)$, as it is a direct summand of $(C^* \otimes_k C)(P) \cong  C(P)^* \otimes_k C(P).$ In particular, $d_i(P)\colon D_i(P)\to D_{i-1}(P)$ is injective for all $P \in s_p(G)$, so $d_i\colon D_i \to D_{i-1}$ is split injective by \cite[Proposition 5.8.11]{L181}. Therefore $D$ has a contractible direct summand, a contradiction. Thus, $C^* \otimes_k C \simeq M[0]$ for some $kG$-module $M$. The K\"unneth formula implies there exists a unique $i\in \Z$ for which $H_i(C)\neq 0$. We conclude $C$ is an endosplit $p$-permutation resolution.
    \end{proof}

    \begin{prop}{\cite[Theorem 7.11.2]{L182}\label{prop:summandsofepprs}}
        Let $M$ be a $kG$-module with endosplit $p$-permutation resolution resolution $C_M$. For any subgroup $H\leq G$, the complex $\End_{kH}(C)$ has homology concentrated in degree zero, isomorphic to $\End_{kH}(M)$, and there is a $k$-algebra isomorphism \[\rho_H\colon \End_{K^b({}_{kH}\catmod)}(C) \cong \End_{kH}(M)\] satisfying $\Res^G_H \circ \rho_G = \rho_H \circ \Res^G_H$ and $\tr^G_H \circ \rho_H = \rho_G \circ \tr^G_H$. In particular, $\rho_G$ induces a vertex and multiplicity preserving bijection between the sets of isomorphism classes of indecomposable direct summands of $M$ and noncontractible indecomposable direct summands of the complex $C$.

        In particular, if $M,N$ are $kG$-modules with $M \oplus N$ having an endosplit $p$-permutation resolution $C_{M \oplus N}$, then $C_{M\oplus N} = C_M \oplus C_N$ with $C_M,C_N$ endosplit $p$-permutation resolutions of $M,N$ respectively.
    \end{prop}

    Let $(K,\calO, k)$ be a $p$-modular system. Endosplit $p$-permutation resolutions over $\calO$ are defined in the analogous way. The next proposition of Rickard ensures that endosplit $p$-permutation resolutions over $\calO$ and $k$ are in bijection.

    \begin{prop}{\cite[Proposition 7.1]{R96}}\label{prop:uniquelift}
        Let $G$ be a finite group, and let $M$ be a $kG$-module that has an endosplit $p$-permutation resolution $X_M$. Then there exists a unique lift of $M$ to a $\calO G$-module $\tilde{M}$ with an endosplit $p$-permutation resolution $\tilde{X}_{M}$ lifting $X_M$, i.e. $X_M \cong k \otimes_\calO \tilde{X}_{M}$.
    \end{prop}

    Note that $M$ may not lift uniquely to $\calO$ - rather, the proposition guarantees the existence of a unique lift to $\calO$ with an endosplit $p$-permutation resolution.

    Finally, any endosplit $p$-permutation resolution has a unique direct summand which has no contractible elements. It will be convenient at times to assume that an endosplit $p$-permutation resolution has no contractible summands; this choice is always unique.

    \begin{definition}{\cite[Definition 6.1]{M24b}}
        Let $C$ be a bounded chain complex of $p$-permutation $kG$-modules. $C$ is a \textit{$V$-endosplit-trivial complex} if and only if $C^* \otimes_k C \simeq (k \oplus N)[0]$ for some $V$-projective $kG$-module $N$. Equivalently, $C$ is an endosplit $p$-permutation resolution of a $V$-endotrivial $kG$-module. For shorthand, we refer to $V$-endosplit-trivial complexes as simply \textit{$V$-endotrivial complexes}. $V$-endotrivial complexes over $\calO$ are defined analogously.

        If $V=0$, then we have $C^* \otimes_k C\simeq k[0]$, recovering the definition of an \textit{endotrivial complex}, see \cite[Definition 3.1]{M24a}. These are the invertible objects of the tensor-triangulated category $K^b({}_{kG}\triv)$, and are completely classified in \cite{M24c}.
    \end{definition}

    \begin{observation}
        If $C^* \otimes_k C \simeq (k \oplus N)[0]$, $N$ is necessarily $p$-permutation. Therefore, when considering the class of $V$-endotrivial complexes, $V$ can be replaced with a suitable $p$-permutation $kG$-module $W$.
    \end{observation}

    \begin{prop}{\cite[Propositions 6.6, 6.10]{M24b}}
        Let $C$ be a $V$-endotrivial complex.
        \begin{enumerate}
            \item $C \simeq C_0 \oplus D$, where $C_0$ is an indecomposable $V$-endotrivial complex and $D$ is a direct sum of $V$-projective chain complexes.
            \item $C_0$ has vertex set $\Syl_p(G)$.
        \end{enumerate}
    \end{prop}
    \begin{proof}
        For (a), let $C = C_0 \oplus C_1$, then \[(C_0^*\otimes_k C_0) \oplus (C_1^* \otimes_k C_1) \oplus (C_0^* \otimes_k C_1) \oplus (C_1^* \otimes_k C_0)  \cong (C_0 \oplus C_1)^* \otimes_k (C_0 \oplus C_1) \simeq (k \oplus N)[0].\] Neither of $(C_0^* \otimes_k C_1), (C_1^* \otimes_k C_0)$ can contain $k[0]$ as a direct summand since they are dual, therefore one of the other two summands contains $k[0]$ as a direct summand. Without loss of generality, suppose that it is $C_0^* \otimes_k C_0$. Then $C_0$ is a $V$-endotrivial chain complex and the direct sum of the three other terms in the direct sum is homotopy equivalent to $N'[0]$ for some $V$-projective $kG$-module $N'$. Therefore, all three other terms in the direct sum are $V$-projective as chain complexes.

        Now, consider \[C_1 \otimes_k (k \oplus N'')[0] \cong C_1 \otimes_k (C_0^* \otimes_k C_0) \cong (C_1 \otimes_k C_0^*) \otimes_k C_0.\] The rightmost term in this chain of isomorphisms is $V$-projective since $C_1 \otimes_k C_0^*$ is. Therefore, all direct summands are as well, so in particular, $C_1$ is $V$-projective, as desired.

        For (b), Let $S \in \Syl_p(G)$. Since $(C^* \otimes_k C)(S)\cong k \neq 0$, $C$ contains at least one degree with Sylow vertices. From the decomposition $C \simeq C_0 \oplus D$ shown in (a), $D(S)=0$ since $V$ is absolutely $p$-divisible. Therefore $C_0$ necessarily has vertex $S$, as desired.
    \end{proof}

    \begin{construction}
        We impose an equivalence condition $\sim_V$ on the class of $V$-endotrivial chain complexes of $kG$-modules (analogous to the module-theoretic case) as follows: if $C$ and $D$ are $V$-endotrivial chain complexes, write $C\sim_V D$ if and only if $C$ and $D$ have isomorphic indecomposable $V$-endotrivial summands. Let $\calE_k^V(G)$ denote the resulting set of equivalence classes. $\calE_k^V(G)$ is an abelian group with addition defined by $[C]+[D] := [C\otimes_k D]$. We have $-[C] = [C^*]$ and each equivalence class has a unique indecomposable representative.
    \end{construction}

    \begin{theorem}
        Let $C$ be a bounded chain complex of $p$-permutation $kG$-modules. The following are equivalent.
        \begin{enumerate}
            \item $C$ is a $V$-endotrivial complex.
            \item $C$ is an endosplit $p$-permutation resolution of a $V$-endotrivial $kG$-module.
            \item For all $P \in s_p(G)$, there exists a unique integer $h_C(P) \in \Z$ such that $H_{h_C(P)}(C(P)) \neq 0$, $H_i(C(P)) = 0$ for $i \neq h_C(P)$, and in addition, if $V(P) = 0$ then $\dim_k H_{h_C(P)}(C(P))=1$.
        \end{enumerate}
    \end{theorem}
    \begin{proof}
        The equivalence between (a) and (b) was stated in the previous definition. Throughout this proof, we assume $V$ is $p$-permutation.

        First, suppose $C$ is a $V$-endotrivial complex. Since $C^* \otimes_k C \simeq (k \oplus N)[0]$ for some $V$-projective $kG$-module $N$, \[C(P)^* \otimes_k C(P) \cong (C^*\otimes_k C)(P) \cong  (k \oplus N(P))[0].\] It follows by the K\"unneth formula that $C(P)$ also has homology concentrated in a single degree, $h_C(P)$, and if $V(P) = 0$, then since $N$ is $V$-projective with $V$ absolutely $p$-divisible, $N(P) = 0$. This shows (c).

        Now, suppose (c). By Proposition \ref{prop:equivendosplit}, $C$ is an endosplit $p$-permutation resolution for $M := H_{h_C(1)}(C)$. It suffices to show $M$ is $V$-endotrivial. Since \\ $\dim_k H_{h_C(P)}(C(P))=1$ for all $P \not\in \calX_V$, $C(P)^* \otimes_k C(P) \cong (C^* \otimes_k C)(P) \simeq k[0]$. On the other hand, since $C$ is an endosplit $p$-permutation resolution, $C^* \otimes_k C \simeq (M^* \otimes_k M)[0]$. Thus $M^* \otimes_k M \cong L$, where $L$ is a $p$-permutation $kG$-module for which $L(P) \cong k$ for any $P \not\in \calX_V$.

        Since $V$ is $p$-permutation, $V$-projectivity is equivalent to $\calX_V$-projectivity, where $\calX_V$ is the set of $p$-subgroups for which $V(P) \neq 0$. A $p$-permutation module $N$ is $\calX_V$-projective if and only if $N(P) = 0$ for all $P \not\in\calX_V$, since $\calX_V$ is downward closed as a poset. Therefore, $L = L_0 \oplus L_1$, where $L_0$ is a $p$-permutation $kG$-module with vertex $S$ and $L_1$ is $V$-projective. Since $L_0$ has vertex $S$ and $L_0(S) \cong k$, $L_0 \cong k$ since $L_0(S)$ is the Green correspondent of $L_0$, see \cite[Proposition 5.10.5]{L181}. We conclude that $M$ is $V$-endotrivial, as desired.
    \end{proof}

    The next construction is fundamental in our study of (relatively) endotrivial complexes.

    \begin{construction}{\cite[Definition 9.6]{M24b}}
        If $C$ is a $V$-endotrivial complex, the values $h_C(P)$ in the previous proposition define an integer-valued superclass function on the $p$-subgroups of $G$, $h_C\colon s_p(G) \to \Z$. We denote the set of such superclass functions by $CF(G,p)$, call the functions $h_C$ the \textit{$h$-marks of $C$}, and say $h_C(P)$ is \textit{the $h$-mark of $C$ at $P$}.

        It is straightforward to verify that the map $h\colon \calE_k^V(G) \to CF(G,p)$ is a well-defined group homomorphism. We call $h$ the \textit{$h$-mark homomorphism}, as the values $h_C(P)$ are analogous to the marks of a $G$-set. That is, if $X$ is a $G$-set, then $X$ is uniquely determined by its \textit{marks} $|X^H|$ where $H$ runs over all conjugacy class representatives of subgroups of $G$; in fact, the assignment $X \mapsto (H \mapsto |X^H|)$ induces an injective ring homomorphism $B(G) \to CF(G)$ (see e.g., \cite[Theorem 2.4.5]{Bou10}), where $B(G)$ denotes the \emph{Burnside ring} of $G$, the split Grothendieck ring of ${}_G\mathbf{set}$.
    \end{construction}

    \begin{notation}
        Let $T_{k,V}(G,S) \leq T_{k,V}(G)$ denote the subgroup of $T_{k,V}(G)$ consisting of $p$-permutation $V$-endotrivial $kG$-modules. Equivalently, if $S\in\Syl_p(G)$, $T_{k,V}(G,S)$ is the kernel of the restriction map $\Res^G_S\colon T_{k,V}(G) \to T_{k,V}(S)$. The group $T_{k,V}(G,S)$ is finite.
    \end{notation}

    \begin{theorem}{\cite[Theorem 9.7]{M24b}, \cite[Theorem 2.12]{M24c}}\label{thm:hmarkhom}
        The homomorphism $h\colon \calE_k^V(G) \to CF(G,p)$ has kernel the torsion subgroup of $\calE_k^V(G)$, which is identified with $T_{k,V}(G,S)$ via the inclusion $[M] \mapsto [M[0]]$. In particular, if $G$ is a $p$-group, $h$ is injective. If $V = V(\calF_G)$, $h$ is surjective.
    \end{theorem}
    \begin{proof}
        The kernel of $h$ consists of (equivalence classes of) $V$-endotrivial chain complexes $C$ for which $C(P)$ has nonzero homology only at degree 0 for any $p$-subgroup $P$. Suppose for contradiction that there exists an indecomposable $V$-endotrivial chain complex $C$ with $[C] \in \ker h$ such that $C \not\cong M[0]$ for some $kG$-module $M$. Without loss of generality suppose that there exists $i > 0$ for which $C_i \neq 0$, and $C_j = 0$ for all $j > i$. Then for all $p$-subgroups $P$, $d_i(P)\colon C_i(P) \to C_{i-1}(P)$ is injective, hence $d_i$ is split injective, a contradiction to $C$ being indecomposable. Thus if $C$ is indecomposable and $[C] \in \ker h$, $C \cong M[0]$ for a necessarily $p$-permutation $V$-endotrivial $kG$-module $M$. Thus, $\ker h \cong T_{k,V}(G,S)$. Since $T_{k,V}(G,S)$ is finite, $\ker h$ is torsion, and any $[C] \not\in \ker h$ clearly cannot be torsion since for any $n \in \N$, $n\cdot h_C$ outputs a nonzero value.

        If $G$ is a $p$-group, then $T_{k,V}(G,S) = \{[k]\}$ since the only transitive permutation $kG$-module with vertex $S = G$ is the trivial one. If $V = V(\calF_G)$, then one may verify that the image of the set \[\{k[G/P] \to k\}_{P \in [s_p(G)\setminus \Syl_p(G)]/G}\cup \{k[1]\}\] under $h$ forms a diagonal basis of $CF(G,p) = CF(G)$. Explicitly, setting $C_P:=k[G/P]\to k$, we have \[h_{C_P}(H) = \begin{cases} 1 & H \leq P\\ 0& H \not\leq P\end{cases}\] and a M\"obius inversion argument shows this is a $\Z$-basis for $CF_b(G)$. Hence $h$ is surjective, and in fact, the above set is a basis for $\calE_k^{V(\calF_G)}(G)$.
    \end{proof}

    \begin{remark}
        It is clear that every $kG$-module that has an endosplit $p$-permutation resolution is endo-$p$-permutation, but, in general, not every endo-$p$-permutation module has an endosplit $p$-permutation resolution. Mazza proved in \cite{Ma03} when $p$ is odd and $G$ is a $p$-group, every endopermutation $kG$-module (hence endo-$p$-permutation $kG$-module) has an endosplit $p$-permutation resolution. On the other hand, if $G$ is a $2$-group, not every endopermutation $kG$-module has an endosplit $p$-permutation resolution; for instance, the unique endotrivial $kQ_8$-module that does not descend to $\F_2$ has no endosplit $p$-permutation resolution.
    \end{remark}

    \begin{definition}\label{not:tensorsofrelsyz}
        \begin{enumerate}
            \item We define the \textit{partial $h$-marks} for any endosplit \\ $p$-permutation resolution $C$ as follows. Let $\calX_C$ denote the set of all $p$-subgroups of $G$ for which $C(P)\not\simeq 0$. It is easy to see that $\calX_C$ forms a downward closed sub-poset of $s_p(G)$ which is closed under $G$-conjugation. For each $P \in \calX_C$, write $h_C(P)$ for the unique integer for which $H_{h_C(P)}(C(P)) \neq 0$. This defines a superclass function on $\calX_C$, $h_C\colon \calX_C\to \Z$.
            \item Say a $kG$-module $M$ is \textit{a module tensor of relative syzygies} if there exist $p$-subgroups $P_1,\dots, P_i$ and $Q_1,\dots, Q_j$, all of which are non-Sylow, such that \[M \cong \Delta(G/P_1)\otimes_k\cdots \otimes_k \Delta(G/P_i)\otimes_k \Delta(G/Q_1)\inv \otimes_k \cdots \otimes_k \Delta(G/Q_j)\inv.\]
            If so, $M$ is strongly capped endo-$p$-permutation, and has an endosplit $p$-permutation resolution
            \begin{align*}
                C_M := &(k[G/P_1]\to k)\otimes_k \cdots \otimes_k (k[G/P_i]\to k)\\
                &\otimes_k (k \to k[G/Q_1])\otimes_k \cdots \otimes_k (k \to k[G/Q_j]).
            \end{align*}
            We have that $C_M$ is a $V(\calF_G)$-endotrivial chain complex; we say a complex homotopy equivalent to a shift of a complex in the above form is \textit{a chain tensor of relative syzygies.} Given any chain tensor of relative syzygies $C$, there exists a unique corresponding module tensor of relative syzygies, obtained by taking the nonzero homology $H_{h_C(1)}(C)$. Conversely, it is clear from the definition that given any module $M$ tensor of relative syzygies, there exists some chain tensor of relative syzygies $C_M$ resolving $M$, however $C_M$ need not be unique (even up to contractibles).

            From the proof of Theorem \ref{thm:hmarkhom}, $\calE_k^{V(\calF_G)}(G)$ has a $\Z$-basis given by (images of) tensors of relative syzygies (including the chain complex $k[1]$). Therefore, the subgroup $\calE_k^{V(\calF_G)}(G)^\Omega$ of $\calE_k^{V(\calF_G)}(G)$ consisting of equivalences classes of tensors of relative syzygies is isomorphic to $CF(G,p)$.
        \end{enumerate}

    \end{definition}

    \begin{remark}
        We are about to characterize all $kG$-modules which have endosplit $p$-permutation resolutions, as well as all endosplit $p$-permutation resolutions themselves. This may be seen as the closest one may get to a classification theorem in the spirit of Dade, as there is no group structure that parameterizes all endo-$p$-permutation $kG$-modules (which have endosplit $p$-permutation resolutions) or all endosplit $p$-permutation resolutions, due to a lack of invertibility of some objects, in a sense.

        Similar to how capped endopermutation modules can be thought of as the invertible endopermutation modules, and thus can be given a group structure, the $V(\calF_G)$-endotrivial complexes can be thought of as the invertible, or ``capped,'' endosplit $p$-permutation resolutions. The following theorem morally states that every endo-$p$-permutation module (resp. endosplit $p$-permutation resolution) can be built from invertible endo-$p$-permutation modules (resp. endosplit-$p$-permutation resolutions), similar to how all endopermutation modules for $p$-groups arise from capped endopermutation modules (see \cite[Theorem 6.6]{D78a}).
    \end{remark}

    \begin{theorem}\label{thm:charofendosplit}
        Let $k$ be any field of characteristic $p > 0$, and let $M$ be an endo-$p$-permutation $kG$-module. $M$ has an endosplit $p$-permutation resolution $C_M$ if and only if there exists a $V(\calF_G)$-endotrivial $kG$-module $E$ which is a tensor of relative syzygies and a permutation module $X$ such that $M$ is a direct summand of $E\otimes_k X$. If this occurs and if $C_M$ has no contractible direct summands, then $C_M$ is a direct summand of a chain complex $C_E \otimes_k X[0]$, where $C_E$ is a chain tensor of relative syzygies which resolves $E$.

        In particular, every endosplit $p$-permutation resolution is a direct summand of a complex $C_E \otimes_k X[0]$, where $C_E$ is a chain tensor of relative syzygies and $X$ is a permutation $kG$-module with $p$-subgroup stabilizers.
    \end{theorem}

    \begin{proof}
        The reverse direction of the first statement is straightforward: since both $E$ and $X$ have endosplit $p$-permutation resolutions, $E\otimes_k X$ does as well, and by Proposition \ref{prop:summandsofepprs}, $M$ has an endosplit $p$-permutation resolution.

        Now, suppose $M$ has an endosplit $p$-permutation resolution $C_M$. Then since $h\colon \calE_k^{V(\calF_G)} \to CF(G,p)$ is surjective, there exists a $V(\calF_G)$-endotrivial complex $C_E$ for which $h_{C_E}(P) = h_{C_M}(P)$ for all $P \in s_p(G)$ for which $C_M(P)$ is not contractible, and $h_{C_E}(Q) = 0$ for all $Q \in s_p(G)$ for which $C_M(Q)$ is contractible. Moreover, we may take $C_E$ to be a chain tensor of relative syzygies. Since $C_E$ is also an endosplit $p$-permutation resolution, the module which it resolves, $E$, is a module tensor of relative syzygies. Now, $C_E^* \otimes_k C_M$ is an endosplit $p$-permutation resolution for $E^*\otimes_k M$ which has partial $h$-marks entirely zero. It follows by inductive removal of contractible summands (similar to the proof of Proposition \ref{prop:equivendosplit}) that $C_E^* \otimes_k C_M \simeq (E^* \otimes_k M)[0]$. Set $X := E^* \otimes_k M$, then $X$ is a $p$-permutation module since it is (as a chain complex) a direct summand of a chain complex of $p$-permutation modules.

        Now, we have \[C_M \oplus (N[0] \otimes_k C_M) \cong (k \oplus N)[0]\otimes_k C_M\simeq C_E \otimes_k C_E^* \otimes_k C_M \simeq  C_E \otimes_k X[0] ,\] where $N$ is some $V(\calF_G)$-projective $p$-permutation $kG$-module. Extracting homology, we find \[M \oplus N \otimes_k M \cong H_{h_{C_M}(1)}(C_E \otimes_k C_E^* \otimes_k C_M)= E \otimes_k E^* \otimes_k M = E\otimes_k X,\] where the leftmost isomorphism follows from the previous line. Therefore, $M \cong H_{h_{C_M}(1)}(C_M)$ is a direct summand of $E \otimes_k X$. Finally, if $X$ is not permutation, there exists another $p$-permutation $X'$ for which $X \oplus X'$ is permutation (in fact, $X'$ can be chosen such that $X \oplus X'$ has stabilizers only $p$-subgroups of $G$, since every trivial source $kG$-module with vertex $P$ is a direct summand of $k[G/P]$), and replacing $X$ with $X \oplus X'$ completes the proof of the first statement.

        Now, suppose that $C_M$ has no contractible direct summands. It follows from the homotopy equivalence $C_M \oplus (N[0] \otimes_k C_M)\simeq C_E \otimes_k C_E^* \otimes_k C_M$ that $C_M$ is a direct summand of $C_E \otimes_k C_E^* \otimes_k C_M$. Now, we have \[C_E \otimes_k C_E^* \otimes_k C_M \cong C_E \otimes_k (X[0] \oplus K) \cong C_E \otimes_k X[0] \oplus C_E \otimes_k K,\] where $K$ is a contractible complex. Note that $C_E \otimes_k K$ is also contractible. Combining these statements, $C_M$ is a direct summand of $C_E \otimes_k X[0] \oplus C_E \otimes_k K$, but since $C_M$ was assumed to have no contractible direct summands, we conclude $C_M$ is a direct summand of $C_E \otimes_k X[0]$.

        The final statement follows since we assumed in the proof that $C_M$ was an arbitrary endosplit $p$-permutation resolution, and from the observation that the permutation module $X$ can be chosen to have $p$-subgroup stabilizers.
    \end{proof}

    In particular, every endosplit $p$-permutation resolution is a direct summand of some chain complex which comes from the base field $\F_p$. If one wishes to construct an endosplit $p$-permutation resolution of a module, it is not necessary to enlarge the base field.

    \begin{corollary}
        Let $k$ be a finite field, and $M$ be an endo-$p$-permutation $kG$-module. Set $k':= \overline{k}$. The module $k' \otimes_k M$ has an endosplit $p$-permutation over $k'$ if and only if $M$ has an endosplit $p$-permutation resolution over $k$.
    \end{corollary}
    \begin{proof}
        The reverse implication is clear. Suppose $k' \otimes_k M$ has an endosplit $p$-permutation resolution. By Theorem \ref{thm:charofendosplit}, $k'\otimes_k M$ is a direct summand of $E\otimes_k X$, where $E$ is a module tensor of relative syzygies and $X$ is a permutation module, both over $k'$. By construction, both $E$ and $X$ have corresponding $kG$-modules $\tilde{E}$ and $\tilde{X}$ for which $k' \otimes_k \tilde{E} \cong E$ and $k' \otimes_k \tilde{X} \cong X$. Then \[k' \otimes_k (\tilde{E}\otimes_k \tilde{X}) \cong (k' \otimes_k \tilde{E}) \otimes_k' (k' \otimes_k \tilde{X})\cong E\otimes_k X.\] By Theorem \ref{thm:charofendosplit}, it suffices to show that $M$ is a direct summand of $\tilde{E} \otimes_k \tilde{X}$. However, this follows by the Noether-Deuring theorem since $k' \otimes_k M$ is a direct summand of $E \otimes_{k'} X \cong k'\otimes_k (\tilde{E} \otimes_k \tilde{X})$.

    \end{proof}

    The classification can be restated in a slightly more compact manner, which is reminiscent of Dade's original work on endopermutation modules. The following corollary also gives another moral reason why $V(\calF_G)$-endotrivial complexes may be seen as the ``capped'' endosplit $p$-permutation resolutions.

    \begin{corollary}\label{cor:dadechar}
        Let $C$ be an endosplit $p$-permutation resolution with no contractible direct summands of a $kG$-module. Then $C$ is a direct summand of the endosplit $p$-permutation resolution \\$\bigoplus_{i=1}^l \Ind^G_{P_i} \Res^G_{P_i} Y$ for some chain tensor of relative syzygies $Y$ and family of $p$-subgroups $\{P_i\}_{i=1}^l$ of $G$.
    \end{corollary}
    \begin{proof}
        We see from Theorem \ref{thm:charofendosplit} that there exists a chain tensor of relative syzygies $Y$ and a permutation $kG$-module $X$ with $p$-subgroup stabilizers such that $C$ is a direct summand of $X[0] \otimes_k Y$. Frobenius reciprocity asserts that \[\bigoplus_{i=1}^l\Ind^G_{P_i}((\Res^G_{P_i} Y) \otimes_k k[0]) \cong \bigoplus_{i=1}^l Y  \otimes_k (\Ind^G_{P_i} k)[0] \cong Y   \otimes_k \left(\bigoplus_{i=1}^{l}\Ind^G_{P_i} k\right)[0],\] and the result follows.
    \end{proof}

    \begin{remark}
        Dade described how all indecomposable endopermutation $kG$-modules arise from capped endopermutation $kG$-modules in a similar manner for $p$-groups. If $G$ is a $p$-group, \cite[Theorem 6.6]{D78a} implies that if $M$ is an indecomposable endopermutation $kG$-module with vertex $Q$, then $M$ is a direct summand of $\Ind^G_Q \Res^G_Q M$. Moreover, \cite[Theorem 6.10]{D78a} gives a way to decompose any capped endopermutation $kG$-module in terms of its cap. These facts along with the previous corollary implies the following.
    \end{remark}

    \begin{corollary}{\cite{Ma03}}
        Let $G$ be a $p$-group. If $D_k(G) = D^\Omega(G)$, then every indecomposable endopermutation $kG$-module and every capped endopermutation $kG$-module has an endosplit $p$-permutation resolution.
    \end{corollary}

    This is a rediscovery of Mazza's theorem, stating when $p$ is odd, every capped endopermutation $kG$-module has an endosplit $p$-permutation resolution. The question of whether, if $D_k(G) = D^\Omega(G)$, any endopermutation $kG$-module has an endosplit $p$-permutation resolution, is a bit more technical, since the decomposition of such modules is more cumbersome, see \cite[Proposition 6.13]{D78a}.

    \section{Galois descent of endosplit $p$-permutation resolutions}

    In this section, we consider the question of Galois descent of endosplit\\ $p$-permutation resolutions. First, we reduce the question of Galois descent of relatively endotrivial complexes to Galois descent of the module that they resolve. Using this, and the classification of endosplit $p$-permutation resolutions, we extend this reduction to all endosplit $p$-permutation resolutions.

    \begin{notation}
        We adapt the notation conventions from \cite{KL18}. For arbitrary extensions of commutative rings $R \subseteq R'$, if $A$ is an $\calO$-algebra, one can form the $R'$-algebra $A' := R' \otimes_R A$. Similarly, given any $A$-module $M$, one can form an $A'$-module $ R' \otimes_R M$. In general, not every $A'$-module arises in this way. Given a ring automorphism $\sigma$ of $R'$ which restricts to the identity map on $R$ and an $A'$-module $M'$, we denote by ${}^\sigma M'$ the $A'$-module which is equal to $M'$ over the subalgebra $1 \otimes A$ of $A'$, and such that $\lambda \otimes a \in A'$ acts on $M'$ as $\sigma\inv(\lambda) \otimes a$. Note that if $f\colon M' \to N'$ is an $A'$-module homomorphism, then $f$ is also an $A'$-module homomorphism $f\colon {}^\sigma M' \to {}^\sigma N'$. Thus, the Galois twist ${}^\sigma(-)$ induces an $R$-linear (but in general not $R'$-linear), exact autoequivalence of ${}_{A'}\catmod$.

        We have functors \[-_{R}\colon {}_{A'}\catmod \to {}_{A}\catmod \text{ and } R' \otimes_R -\colon {}_{A}\catmod \to {}_{A'}\catmod,\] restriction and extension of scalars, respectively. These are both $R$-linear exact functors. Moreover, $R'\otimes_R -$ is left adjoint to $ -_R$. Moreover, these functors induce functors over the bounded chain complex categories $Ch^b({}_{A}\catmod)$ and $Ch^b({}_{A'}\catmod)$ and the same adjunction holds.

        We say an $A'$-module $M$ (resp. chain complex of $A'$-modules $C$) \textit{descends to $A$} or \textit{descends to $R$} if there exists a $A$-module $\tilde{M}$ (resp. chain complex of $A$-modules $\tilde{C}$) such that $R' \otimes_R \tilde{M} \cong M$ (resp. $R' \otimes_R \tilde{C} \cong C$). We may also say $M$ \textit{is realized in $R$}. Another way of viewing this (see \cite[Remark after 7.15]{CR81}) is that there exists a matrix representation of $M$ over $R'$ with coefficients in $R$. We have the same notion for homomorphisms of $R'G$-modules, by considering them as intertwiners of matrix representations.

        Let $(K, \calO, k) \subseteq (K', \calO', k')$ be an extension of $p$-modular systems, i.e. a pair of $p$-modular systems such that $\calO \subseteq \calO'$ and $J(\calO) \subseteq J(\calO')$. Furthermore, we assume that $k$ and $k'$ are finite fields, and that $\calO$ and $\calO'$ are absolutely unramified, i.e. $J(\calO) = p\calO$ and $J(\calO') = p\calO'$. This is always achievable by setting $\calO = W(k)$ and $\calO' = W(k')$, the Witt vectors of $k$ and $k'$. Set $d = [k' : k]$. Then $\calO'$ is free of rank $d$ as an $\calO$-module. Let $\sigma\colon k' \to k'$ be a generator of $\Gal(k'/k)$. Then there exists a unique ring homomorphism $\sigma\colon \calO' \to \calO'$ lifting $\sigma\colon k' \to k'$, since $\calO'$ is absolutely unramified, and the composite $(\calO' \otimes_\calO - )_\calO$ is naturally isomorphic to the functor $(-)^{\oplus d}$ since $\calO'$ is free over $\calO$ of rank $d$. As a result, since the Krull-Schmidt theorem for $\calO G$-modules and complexes holds, Galois descent is unique up to isomorphism. That is, if $M_1, M_2$ are $\calO G$-modules for which $\calO' \otimes_\calO M_1 \cong \calO' \otimes_\calO M_2$, then $M_1 \cong M_2$, and similarly for chain complexes. This also holds over the residue field $k$.

    \end{notation}

    Our main tool will be the following theorem, which describes exactly when Galois descent occurs for chain complexes for algebras over absolutely unramified complete discrete valuation rings.

    \begin{theorem}{\cite[Theorem 6]{M24d}}\label{thm:galoisstablelift}
        Let $(K, \calO, k) \subseteq (K', \calO', k')$ be an extension of $p$-modular systems with $\calO, \calO'$ absolutely unramified, and let $A$ be a finitely generated $\calO$-free $\calO$-algebra. Set $A' := \calO' \otimes_\calO A$.
        \begin{enumerate}
            \item Suppose $C \in Ch^b({}_{A'}\catmod)$ indecomposable satisfies ${}^\sigma C \cong C$ for all $\sigma \in\Gal(k'/k)$, where we regard $\sigma$ as the unique ring homomorphism of $\calO'$ lifting $\sigma \in\Gal(k'/k)$. Then there exists an indecomposable chain complex $\tilde{C} \in Ch^b({}_{A}\catmod)$ such that $\calO' \otimes_\calO \tilde{C} \cong C$. Moreover, $\tilde{C}$ is unique up to isomorphism.
            \item Conversely, let $\tilde{C} \in Ch^b({}_{A}\catmod)$ and define $C := \calO' \otimes_\calO \tilde{C}$. Then $C$ satisfies ${}^\sigma C \cong C$ for all $\sigma \in\Gal(k'/k)$.
        \end{enumerate}
    \end{theorem}

    \begin{lemma}
        Let $\sigma \in \Gal(k'/k)$, $V\in {}_{k'G}\catmod$, $M$ an indecomposable relatively $V$-endotrivial $k'G$-module, and $C$ an indecomposable endosplit $p$-permutation for $M$, i.e. $C$ is a $V$-endotrivial complex. We have ${}^\sigma M \cong M$ if and only if ${}^\sigma C \cong C$.
    \end{lemma}
    \begin{proof}
        First, assume ${}^\sigma M \cong M$. Since $C$ is an endosplit $p$-permutation resolution for $M$ and ${}^\sigma(-)$ is exact, ${}^\sigma C$ is an endosplit $p$-permutation resolution for ${}^\sigma M$. Therefore, $C^* \otimes_k {}^\sigma C$ is an endosplit $p$-permutation resolution for $M^* \otimes_k {}^\sigma M \cong M^* \otimes_k M \cong k \oplus N$ for some $V$-projective $N$.

        Notice that $h_C(P) = h_{{}^\sigma C}(P)$ for all $P \in s_p(G)$. Therefore, $C^* \otimes_k {}^\sigma C$ has $h$-marks entirely 0 and homology ${}^\sigma M \otimes_k M^* \cong k \oplus N$. Therefore, $C^* \otimes_k {}^\sigma C \simeq (k \oplus N)[0]$, so $[C] = [{}^\sigma C] \in \calE_k^V(G)$. However, both $C$ and ${}^\sigma C$ are indecomposable, and every equivalence class in $\calE_k^V(G)$ has a unique indecomposable representative, thus $C \cong {}^\sigma C$. The reverse implication is clear since ${}^\sigma(-)$ is exact.
    \end{proof}

    \begin{corollary}
        Let $k'/k$ be an extension of finite fields and $V\in {}_{k'G}\catmod$. Suppose $C$ is an indecomposable endosplit $p$-permutation for a $V$-endotrivial module $M \in {}_{k'G}\catmod$, i.e. $C$ is $V$-endotrivial. Then $M$ descends to $k$ if and only if $C$ does.
    \end{corollary}
    \begin{proof}
        Fix $(K, \calO, k) \subseteq (K', \calO', k')$ an extension of $p$-modular systems with $\calO, \calO'$ absolutely unramified (for instance, $\calO := W(k), \calO' := W(k')$). The reverse direction of the if and only if statement is clear. Suppose $M$ descends to a $kG$-module $\tilde{M}$. Then ${}^\sigma M \cong M$ for all $\sigma \in \Gal(k'/k)$. By the previous lemma, ${}^\sigma C \cong C$ for all $\sigma \in \Gal(k'/k)$ as well. Now, there exists a unique lift $\hat{C}$, a chain complex of $\calO' G$-modules, such that $k' \otimes_{\calO'} \hat{C} \cong C$. It follows that ${}^\sigma \hat{C} \cong \hat{C}$ for all $\sigma \in \Gal(k'/k)$ by unique lifting of $\sigma$. By Theorem \ref{thm:galoisstablelift}, there exists a chain complex of $p$-permutation $\calO G$-modules $\tilde{C}$ for which $\calO'\otimes_\calO \tilde{C} \cong \hat{C}$. Then $k \otimes_\calO \tilde{C}$. All together, we have \[k' \otimes_k (k \otimes_\calO \tilde{C}) \cong k' \otimes_\calO \tilde{C}\cong k' \otimes_{\calO'} (\calO' \otimes_\calO \tilde{C}) \cong C,\] thus $C$ descends to $k$, as desired.

    \end{proof}

    \begin{remark}
        If the previous theorem holds, then if $M$ descends to the $kG$-module $\tilde{M}$ and $C$ descends to a chain complex of $p$-permutation $kG$-modules $\tilde{C}$, then $\tilde{C}$ is an endosplit $p$-permutation resolution for $\tilde{M}$, since extension of scalars is exact.     This result also holds over complete discrete valuation rings, even if we make no assumptions about absolute unramification, due to unique lifting of endosplit $p$-permutation resolutions.
    \end{remark}

    Because endotrivial complexes are particular examples of relatively endotrivial complexes, we can discern which fields they descend to rather easily. In particular, up to a twist by a degree one character, endotrivial complexes descend to $\F_p$. We denote the group of endotrivial complexes by $\calE_k(G)$, and $\calT\calE_k(G)$ the \textit{homology-normalized} subgroup of $\calE_k(G)$, which consists of equivalence classes of endotrivial complexes with nonzero homology isomorphic to the trivial $kG$-module.

    \begin{corollary}
        Let $C$ be an indecomposable endotrivial complex of $kG$-modules with $H_{h_C(1)}(C) \cong k_\omega$ for some $\omega \in \Hom(G, \F_p^\times).$ Then there exists an endotrivial complex $\tilde{C} \in Ch^b({}_{\F_pG}\triv)$ for which $k \otimes_{\F_p} \tilde{C} \cong C$. In particular, \[\calE_k(G) \cong \calT\calE_{\F_p}(G) \times \Hom(G,k^\times)\cong CF_b(G,p) \times \Hom(G,k^\times),\] and $\calT\calE_k(G)$ is independent of the field $k$.
    \end{corollary}
    \begin{proof}
        The first statement follows immediately from the previous results. This asserts that $\calT\calE_k(G) \cong \calT\calE_{\F_p}(G)$. The isomorphism $\calE_k(G) \cong \calT\calE_{\F_p}(G) \times \Hom(G,k^\times)$ follows from \cite[Remark 3.8]{M24a}, and the isomorphism $\calT\calE_k(G) \cong CF_b(G,p)$ follows from \cite[Corollary 6.4]{M24c}.
    \end{proof}

    \begin{remark}
        The fact that all endotrivial complexes live in $\F_p$ up to a twist by a degree one Brauer character is not that surprising, since they are determined almost entirely by the Borel-Smith functions which exist independently of scalars. Since endotrivial complexes are, in a sense, the one-sided analogue of a splendid Rickard equivalence, this may provide a heuristic reason for why one may place faith in the refined abelian defect group conjecture. Additionally, the previous corollary gives another proof of non-surjectivity in general of the Lefschetz homomorphism $\Lambda\colon \calE_k(G) \to O(T(kG))$, which was first considered in \cite[Section 7]{M24a}.
    \end{remark}

    Now we extend our results to all endosplit $p$-permutation resolutions by using the characterization of endosplit $p$-permutation resolutions from Theorem \ref{thm:charofendosplit}.

    \begin{theorem}\label{thm:descentofepprs}
        Let $k'/k$ be an extension of finite fields. Let $C$ be an endosplit $p$-permutation resolution for an endo-$p$-permutation $k'G$-module $M$ which contains no contractible summands. Then $C$ descends to $k$ if and only if $M$ does.
    \end{theorem}
    \begin{proof}
        The forward direction is clear since extension of scalars is exact, therefore a descendant of $C$ is a resolution for a descendant of $M$. Suppose $M$ descends to $k$. By Theorem \ref{thm:charofendosplit}, $C$ is a direct summand of the tensor product $C_E \otimes_{k'} X[0]$ of a chain tensor of relative syzygies $C_E$ and a permutation module $X$, regarded as a chain complex in degree zero. Now, observe $C_E$ and $X$ both in fact descend to $\F_p$ (this is clear from their constructions), therefore they both descend to any intermediate extension $k' \supseteq k \supseteq \F_p$. Hence the same holds for $C_E \otimes_{k'}X[0]$.

        Similarly, $M$ is a direct summand of $E\otimes_{k'} X$ which also descends to $k$, say with multiplicity $m$. Let $Y$ be the unique up to isomorphism $kG$-module for which $k' \otimes_k Y \cong E\otimes_{k'} X$ and let $C_Y$ be its corresponding endosplit $p$-permutation resolution containing no contractible direct summands. Both complexes are tensors of relative syzygies, and we have $k' \otimes_k C_Y \simeq C_E \otimes_{k'} X[0]$.

        Let $\tilde{M}$ be the $kG$-module for which $k' \otimes_k \tilde{M} \cong M$, then $\tilde{M}$ is a direct summand of $Y$, and also occurs with multiplicity $m$.

        Let $Y = Y_1 \oplus \cdots \oplus Y_l$ be a direct sum decomposition of $Y$ into indecomposables. By Proposition \ref{prop:summandsofepprs}, this corresponds to a direct sum decomposition $C_Y = C_{Y_1} \oplus\cdots \oplus C_{Y_l}$ where each $C_{Y_i}$ is an endosplit $p$-permutation resolution for $Y_i$ with no contractible summands. Therefore, \[C_E \otimes_{k'} X[0] \simeq k' \otimes_k ( C_{Y_1} \oplus\cdots \oplus C_{Y_l}) \cong k' \otimes_k C_{Y_1} \oplus \cdots \oplus k' \otimes_k C_{Y_l}.\] Note that in this sum, there must be exactly $m$ endosplit $p$-permutation resolutions of $M$ present (with respect to this decomposition into not necessarily indecomposable parts), since $\tilde{M}$ occurs in $Y$ with multiplicity $m$ by the Noether-Deuring theorem. Since there must be a total of $m$ endosplit $p$-permutation resolutions of $M$ in any decomposition of $C_E\otimes_{k'} X[0]$, by Krull-Schmidt one of these resolutions must be $C$ up to isomorphism, since $C$ is a direct summand of $C_E \otimes_{k'} X[0]$ as well. Thus $C$ also descends to $k$.

    \end{proof}

    The above theorem also holds for complete discrete valuation rings with finite residue field by the unique lifting of endosplit $p$-permutation resolutions.

    \begin{corollary}
        Let $k$ be a finite field and set $k' := \overline{k}$. Suppose $M$ is an endo-$p$-permutation $kG$-module for which $k' \otimes_k M$ has an endosplit $p$-permutation resolution $\hat{C}$ of $p$-permutation $k'G$-modules. Then there exists an endosplit $p$-permutation resolution $C$ of $M$ such that $k' \otimes_k C \cong \hat{C}$.
    \end{corollary}
    \begin{proof}
        This follows since both every component $\hat{C}_i$ or differential $d_i\colon \hat{C}_i \to \hat{C}_{i-1}$ has some finite field $l$ that realizes $\hat{C}_i$ or $d_i$, since $k'$ is the colimit of the diagram $\F_p \hookrightarrow \F_{p^2} \hookrightarrow \dots$, so any matrix representation or intertwiner will have values contained in some finite field. Therefore, there exists a finite field $l$ which realizes $\hat{C}$, i.e. there exists a chain complex $C'$ of $lG$-modules for which $k' \otimes_l C' \cong \hat{C}$. By Theorem \ref{thm:descentofepprs}, $C'$ descends to $k$, as desired.
    \end{proof}

    \section{Constructing splendid Rickard equivalences for blocks of $p$-solvable groups over $\F_p$}

    In this section, we apply the results of the previous section to show that Kessar--Linckelmann's refinement of Brou\'{e}'s abelian defect group conjecture holds for certain blocks of $p$-solvable groups $G$ satisfying $G = O_{p',p,p'}(G)$ with abelian defect groups.

    \begin{conjecture}{\cite{KL18}}
        ``Given a categorical equivalence, say a Morita equivalence or a Rickard equivalence between $\calO' Gb$ and $\calO' H c$ for some complete discrete valuation ring $\calO'$ for $G$ and $H$ finite groups, $b$ and $c$ blocks of $\calO' G$ and $\calO' H$ respectively, and a complete discrete valuation ring $\calO$ contained in $\calO'$, is the equivalence between $\calO' Gb$ and $\calO' H c$ an extension of an equivalence between $\calO Gb$ and $\calO Hc$? ... [This] may be viewed as evidence for a refined version of the Abelian defect group conjecture, namely that for any $p$-modular system $(K, \calO, k)$ and any block $b$ of $\calO G$ with abelian defect group $P$ and Brauer correspondent $c$, there is a splendid Rickard equivalence between $\calO Gb$ and $\calO N_G(H) c$.''
    \end{conjecture}

    We do so by showing that Harris--Linckelmann's construction in \cite{HL00} can be adapted to a $p$-modular system with $k = \F_p$. One may use unique lifting of splendid Rickard equivalences to lift to $\calO = \Z_p$, which then implies the result for any extension of the $p$-modular system $(\Q_p, \Z_p, \F_p)$. In \cite{HL00}, the authors specify a $p$-modular system with $k$ algebraically closed. However, for many of the statements, one may replace the assumption that $k$ is algebraically closed with the assumption that $k$ is perfect with enough roots of unity, which we will use frequently.

    For this section, we will refer to ``blocks'' as the primitive idempotents of $Z(\overline{\F}_pG)$ (note these always exist over some finite field), and refer to block idempotents over a finite field $k$ which may not be primitive upon extending scalars as ``block idempotents of $kG$.'' Note that for these block idempotents, defect groups are defined in the usual way (for instance, a maximal $p$-subgroup $P$ of $G$ for which $\Br_P(b) \neq 0$) and satisfy essentially all the usual properties. Brauer's first main theorem asserts that there is a canonical bijection between the set of block idempotents of $kG$ with defect group $P$ and the set of block idempotents of $kN_G(P)$ with defect group $P$. Given a block idempotent $b$ of $kG$ with defect group $P$, the block idempotent $c$ of $kN_G(P)$ corresponding to $b$ under this bijection is called the \textit{Brauer correspondent.}

    Given a finite group $G$, we denote by $O_p(G)$ (resp. $O_{p'}(G)$) the largest normal $p$-subgroup (resp. $p'$-subgroup) of $G$. $O_{p',p}(G)$ denotes the inverse image in $G$ of $O_p(G/O_{p'}(G))$, and $O_{p',p,p'}(G)$ denotes the inverse image in $G$ of $O_{p'}(G/O_{p',p}(G))$.

    We begin by recalling some facts about $\calO^\times$-groups and twisted algebra constructions, following \cite{HL00}. For this section, we let $(K, \calO, k)$ be a $p$-modular system, with the possible assumption that $K = \calO = k$ unless otherwise specified. We will note what assumptions are needed on our $p$-modular system, as many statements do not require $k$ to be algebraically closed or even to have enough roots of unity.

    \begin{definition}
        An $\calO^\times$-group is a group $\hat{G}$ endowed with a group homomorphism $\calO^\times \to \hat{G}$ whose image lies in $Z(\hat{G})$; we usually denote by $\hat{\lambda}$ the image of $\lambda \in \calO^\times$ in $\hat{G}$, and by $\check{G}$ the \textit{opposite} $\calO^\times$-group, which as an abstract group is equal to $\hat{G}$ but endowed with the homomorphism $\calO^\times \to \check{G}$ sending $\lambda \in \calO^\times$ to $\check{\lambda} = (\hat{\lambda})\inv.$
    \end{definition}

    \begin{construction}
        For our context, $\calO^\times$-groups arise from the action of $G$ on a matrix algebra $S$ over $\calO$ by algebra automorphisms. Denote by ${}^xs$ the action of $x \in G$ on $s \in S$. By the Skolem-Noether theorem, any automorphism of $S$ is inner, and hence for any $x \in G$, there is $s \in S$ such that ${}^xt = sts\inv$ for all $t \in S$. The group $\hat{G}$ of all such pairs $(x,s) \in G\times S^\times$ satisfying ${}^xt = sts\inv$ for all $t \in S$ becomes an $\calO^\times$-group with the group homomorphism sending $\lambda \in \calO^\times$ to $\hat\lambda = (1_G, \lambda 1_S)$. We call $\hat{G}$ \textit{the $\calO^\times$-group defined by the action of $G$ in $S$.} In fact, $\hat{G}$ is a central $\calO^\times$-extension of $G$, as we have a short exact sequence \[1 \to \calO^\times \to \hat{G}\to G \to 1.\]

        One instance in which this occurs is the following: if $b$ is a block of $G$ with trivial defect group and $\calO$ is large enough, the block algebra $\calO Gb$ is a matrix algebra, and therefore admits a $\calO^\times$-group which we denote $\hat{G}_{\calO, b}$. If $\calO \subseteq \calO'$, the inclusion induces a monomorphism of short exact sequences:

        \begin{figure}[H]
            \begin{tikzcd}
                1 \ar[r] & \calO^\times \ar[r] \ar[d, hookrightarrow] & \hat{G}_{\calO, b} \ar[r] \ar[d, hookrightarrow] & G \ar[r] \ar[d, "="]& 1\\
                1 \ar[r] & \calO'^\times \ar[r] & \hat{G}_{\calO', b}  \ar[r]& G \ar[r]& 1
            \end{tikzcd}
        \end{figure}

        Here, the middle inclusion is induced by the inclusion $G \times \calO Gb\to G\times \calO' Gb$.
    \end{construction}

    \begin{prop}{\cite[Corollary 5.3.4]{L181}}\label{prop:hatpsplitscanonically}
        Let $(K, \calO, k)$ be a $p$-modular system with $k$ algebraically closed. If $S$ has $\calO$-rank prime to $p$ and $P$ is a finite $p$-group acting on a matrix algebra $S$, the action of $P$ lifts to a unique group homomorphism $\sigma\colon P \to S^\times$ such that $\det(\sigma(x)) = 1$ for all $x \in P$. In particular, $\hat{P}$ splits canonically as a direct product $\hat{P} \cong P \times \calO^\times$, with the map sending $u \in P$ to $(u, \sigma(u)) \in \hat{P}$ being a section of the canonical surjection $\hat{P} \to P$.
    \end{prop}

    \begin{remark}
        For every $S$ and $P$ satisfying the hypotheses in the previous proposition, there exists a $p$-modular system $(K, \calO, k)$ with $k$ finite for which the previous proposition holds as well. See \cite[Remark 5.3.11]{L181} and \cite[Remark 1.2.19]{L181} for more details.
    \end{remark}

    \begin{definition} \label{def:centextensionalg}
        If $\hat{G}$ is an $\calO^\times$-group, we denote by $\calO_* \hat{G}$ the quotient of the group algebra $\calO \hat{G}$ by the ideal generated by the set of elements $\lambda.\hat{x} - 1_\calO.(\hat{\lambda}\hat{x})$ where $\lambda$ runs over $\calO$ and $\hat{x}$ runs over $\hat{G}$. $\calO_* \hat{G}$ is $\calO$-free of rank $|G|$, since the image in $\calO_\ast \hat{G}$ of any subset of order $|G|$ of $\hat{G}$ mapping onto $G$ under the canonical surjection $\hat{G} \to G$ is an $\calO$-basis of $\calO_* \hat{G}$.

        Let $(K, \calO, k)\subseteq (K', \calO', k')$ be an extension of $p$-modular systems, let $S$ be a matrix algebra over $\calO$, and set $S' := \calO' \otimes_\calO S$. Let $\hat{G}_\calO$ (resp. $\hat{G}_{\calO'}$) denote the $\calO^\times$-group (resp. $(\calO')^\times$-group) defined by the action of $G$ on $S$ (resp. $G$ on $S'$). Then, the obvious inclusion of $\calO'$-algebras $\calO' \otimes_\calO \calO_* \hat{G}_\calO \to \calO_*'\hat{G}_{\calO'}$ is an isomorphism, since it sends the $\calO'$-basis given by representatives of $G$ in $\hat{G}$ to the analogous $\calO'$-basis.
    \end{definition}

    \begin{prop}{\cite[Statements 2.5, 2.6]{HL00}}\label{prop:indcommute}
        If $\hat{G}$ is a central $\calO^\times$-extension of a finite group $G$, there is a unique natural algebra homomorphism $\calO G \to \calO_\ast \hat{G} \otimes_\calO \calO_\ast \check{G}$ mapping $x \in G$ to $\hat{x} \otimes \hat{x}$, where $\hat{x}$ is any element of $\hat{G}$ (and hence also of $\check{C}$) which lifts $x$. Denote by $\Ind^{\hat{G} \times \check{G}}_{\Delta G}\colon {}_{\calO G}\catmod \to {}_{\calO_\ast \hat{G} \otimes_\calO \calO_\ast \check{G}}\catmod$ the obvious induction functor. Note any $\calO_\ast \hat{G} \otimes_\calO \calO_\ast \check{G}$-module is equivalently a $(\calO_\ast \hat{G}, \calO_\ast \hat{G})$-bimodule since $(\calO_\ast \check{G})^o \cong \calO_\ast \hat{G}$.

        Furthermore, for any two $\calO G$-modules $U,V$ there is a unique natural isomorphism of $(\calO_*\hat{G}, \calO_*\hat{G})$-bimodules \[\Ind^{\hat{G} \times \check{G}}_{\Delta G}(U) \otimes_{\calO_\ast \hat{G}}\Ind^{\hat{G} \times \check{G}}_{\Delta G}(V) \cong \Ind^{\hat{G} \times \check{G}}_{\Delta G}(U\otimes_\calO V)\] mapping $((\hat{x},1)\otimes u)\otimes((1,\hat{y})\otimes v)$ to $(\hat{x},\hat{y})\otimes (u,v)$, where $\hat{x},\hat{y} \in \hat{G}, u\in U, v\in V$.
    \end{prop}

    To prove Brou\'{e}'s conjecture holds, Harris--Linckelmann provide a reduction to the case where $G = O_{p',p,p'}(G)$ with blocks of principal type, then use structural properties of groups of this form to construct a splendid Rickard autoequivalence built from a Morita equivalence and a diagonally induced endosplit $p$-permutation resolution. We verify that their construction descends to the block's field of realization. In \cite{HL00}, the construction is performed over $\calO$, however, we will work in positive characteristic instead, then invoke Rickard's lifting theorem, \cite[Theorem 5.2]{R96}.

    We first restate the structural theorem that describes block decompositions in the restricted case. In the original paper, this is stated over $\calO$, but it is straightforward that this also holds over $k$.

    \begin{prop}{\cite[Proposition 3.5]{HL00}}\label{prop:structureforreducedcase}
        Assume $k$ is algebraically closed. Let $G$ be a finite group with abelian Sylow $p$-subgroup $P$ such that $G = O_{p',p,p'}(G)$. Let $b$ be a $G$-stable block of $O_{p'}(G)$. Set $H = N_G(P)$ and $S = kO_{p'}(G)b$. Denote by $\hat{G}$ the $k^\times$-group opposite to that defined by the action of $G$ on $S$, and set $\hat{L} = \hat{G}/O_{p'}(G)$, where we identify $O_{p'}(G)$ with its natural image $\{(x,xb)\}_{x \in O_{p'}(G)}$ in $\hat{G}$.
        \begin{enumerate}
            \item Then $b$ is a block of $G$ having $P$ as a defect group, and there is a unique algebra isomorphism
            \begin{align*}
                kGb &\cong S\otimes_k k_\ast \hat{L}\\
                xb &\mapsto s\otimes \overline{(x,s)}
            \end{align*}
            where $x\in G$, $s\in S^\times$ such that $(x,s)\in \hat{G}$ and $\overline{(x,s)}$ is its image in $\hat{L}$.

            \item The action of $P$ on $S$ lifts to a unique group homomorphism $\sigma\colon P \to S^\times$ such that $\det(\sigma(x)) = 1$ for all $x \in P$, and there is a $p'$-subgroup $E$ of $\Aut(P)$ and a central $k^\times$-extension $\hat{E}$ of $E$ such that we have an isomorphism of $k^\times$-groups
            \[\hat{L} \cong P \rtimes \hat{E}\]
            mapping $\overline{(u, \sigma(u))}$ to $u$ for all $u \in P$. In particular, the isomorphism in (a) maps $ub$ to $\sigma(u)\otimes \overline{(u, \sigma(u))}$ for any $u \in P$.

            \item Let $c$ be the block of $H$ corresponding to $b$. Then $c$ is a block of $O_{p'}(H)$, we have $O_{p'}(H) = O_{p'}(C_G(P))\subseteq O_{p'}(G)$, and setting $T = kO_{p'}(H)c$ and denoting by $\hat{H}$ the $k^\times$-group oppose to that defined by the action of $H$ on $T$, the inclusion $H \leq G$ induces an isomorphism $H/O_{p'}(H)\cong G/O_{p'}(G)$ which lifts to an isomorphism of $k^\times$-groups $\hat{H}/O_{p'}(H)\cong \hat{G}/O_{p'}(G)$; in particular, we again have an algebra isomorphism \[kHc\cong T\otimes_k k_\ast \hat{L},\] which, for any $u \in P$, maps $uc$ to $1_T \otimes \overline{(u, \sigma(u))}$.
        \end{enumerate}
    \end{prop}

    \begin{remark}\label{rem:reduction}
        The assumption that $k$ is algebraically closed is invoked for Proposition \ref{prop:hatpsplitscanonically}, \cite[Lemma 2.7]{HL00} which is a consequence of Dade's splitting theorem of fusion in endopermutation modules, and in the assumption that $S$ is a matrix algebra. One may verify Dade's theorem satisfies the same large enough finite field replacement condition via a similar argument; see the proof of \cite[Theorem 7.9.1]{L182} and \cite[Remark 1.2.19]{L181}. Therefore, there exists a $p$-modular system $(K,\calO, k)$ with $k$ finite (large enough for $G$) for which the previous proposition holds as well.

        In the setup of (c), since $P$ acts trivially on $T$, $T$ is source algebra equivalent (hence splendidly Morita equivalent) to $k$. We claim this equivalence descends to the subfield $\F_p[c]$. Indeed, Wedderburn's theorem asserts $\F_p[c]O_{p'}(H)c \cong M_n(l)$ (as $\F_p[c]$-algebras) for some finite field extension $l/\F_p[c]$, but $\F_p[c]O_{p'}(H)c$ remains indecomposable upon tensoring by any field extension of $\F_p[c]$, and this happens only for $l = \F_p[c]$.

        Therefore, we have a splendid Morita equivalence between $kHc$ and $k_\ast \hat{L}$, since source algebra equivalences induce splendid Morita equivalences (see \cite[Theorem 4.1]{L01}, note that the notion of splendid-ness used in this article is stronger than Rickard's original notion). Therefore, to verify Brou\'{e}'s conjecture, it suffices to construct a splendid Rickard equivalence for $kGb$ and $k_\ast \hat{L}$. To verify the refined abelian defect group conjecture, it suffices to show that the equivalence descends to $\F_p[b] = \F_p[c]$ (see \cite[Lemma 4.2]{BY22}), then apply \cite[Theorem 6.5(b)]{KL18}.
    \end{remark}

    Finally, we recall a classical fact in the classification of endopermutation modules for $p$-groups.

    \begin{lemma}\label{lemma:endopermsforabeliansdescend}
        Let $P$ be an abelian $p$-group and $k$ a finite field. The following hold:
        \begin{enumerate}
            \item $D_k(P) = D^\Omega(P)$, where $D^\Omega(P)$ is the subgroup of $D_k(P)$ generated by equivalence classes of relative syzygies, that is, kernels of augmentation homomorphisms $kX \to k$, where $X$ is a transitive $P$-set.
            \item If $M$ is an indecomposable endopermutation $kP$-module, then $M$ is absolutely indecomposable.
        \end{enumerate}
        In particular, every endopermutation $kP$-module descends to $\F_p$ (or any intermediate field extension), that is, there exists an endopermutation $\F_p P$-module $\tilde{M}$ such that $k \otimes_{\F_p} \tilde{M} \cong M$.
    \end{lemma}
    \begin{proof}
        Statement (a) follows from the classification in \cite[Theorem 9.5]{Bo06} (clarifying the original classification in \cite{D78a}), and statement (b) follows from \cite[Theorem 6.6]{D78a}. For the final statement, see \cite[Remark 7.8.5]{L182}.
    \end{proof}

    \begin{construction}\label{cons:descentofspl}
        We now discuss the construction of the splendid Rickard equivalence between $kGb$ and $kHc$ with our assumptions in place and prove descent. We continue to use the notation of Proposition \ref{prop:structureforreducedcase}. As noted in Remark \ref{rem:reduction}, it suffices to construct a splendid Rickard equivalence between $kGb \cong kO_{p'}(G)b \otimes_k k_\ast \hat{L}$ and $k_\ast\hat{L}$. There are two components of this construction, an $S$-module $V$ (which as $kP$-module is endopermutation) and a Rickard tilting complex (non-splendid) of $(k_\ast \hat{L}, k_\ast\hat{L})$-bimodules $X$. We describe the construction of both of them, and show that they both descend to $\F_p[b]$.

        Set $l := \F_p[b]$; we have an extension of finite fields $k/l/\F_p$. Set $\tilde{S} := l O_{p'}(G)b$; it is clear that $k \otimes_l \tilde{S} \cong S$ as $k$-algebras, but we must still show they are isomorphic as interior $P$-algebras. Abusing notation, we construct the $l^\times$-group $\hat{G}$ by considering the action of $G$ on $\tilde{S}$, set $\hat{T} := \hat{G}/O_{p'}(G)$ as in Proposition \ref{prop:structureforreducedcase}, and we have an isomorphism of $k$-algebras $k \otimes_l l_\ast \hat{T} \cong k_\ast \hat{T}$ as noted in Definition \ref{def:centextensionalg}.

        For all details of the original construction, see \cite[Section 4]{HL00}. Note that the original construction is performed over $\calO$, however, we perform the construction over $k$ for simplicity, and use unique lifting of splendid Rickard equivalences from $k$ to $\calO$.

        \begin{itemize}
            \item By Proposition \ref{prop:structureforreducedcase}(b), there is an isomorphism of interior $P$-algebras $S \cong \End_k(V)$ for some $kP$-module $V$ with determinant $1$. Note $kO_{p'}(G)$ is a permutation $kP$-module, where $P$ acts by conjugation. Therefore, $S = kO_{p'}(G)b$ is also a permutation $kP$-module, so $V$ is an endopermutation $kP$-module. Since we have an algebra isomorphism $S \cong \End_k(V)$, $V$ is an $S$-module. We have $\End_k(V) \cong V \otimes_k V^* \cong S$ as $(S,S)$-bimodules. (Note that in the original construction, $V$ is an $\calO P$-module, and another endopermutation $\calO P$-module $W$ is constructed, however $k\otimes_\calO V = k \otimes_\calO W$, so we omit this step.)

            By Lemma \ref{lemma:endopermsforabeliansdescend}, the $kP$-module $V$ descends to $l$ (in fact, it descends to $\F_p$). Let $\tilde{V}$ denote the endopermutation $l P$-module for which $k \otimes_l \tilde{V}\cong V$. Then the structural homomorphism $P \to \End_k(V)^\times$ has image contained in the $l$-subalgebra $\End_l(\tilde{V})$, therefore under the interior $P$-algebra isomorphism $S \cong \End_k(V)$, the image of $\sigma\colon P \to S^\times$ is contained in the $l$-subalgebra $\tilde{S} = lO_{p'}(G)b$. In particular, we have an isomorphism of interior $P$-algebras $S \cong k \otimes_l \tilde{S}$. Therefore, we have an isomorphism of interior $P$-algebras $\tilde{S} \cong \End_l(\tilde{V})$. It follows again that $\tilde{V}$ is an $\tilde{S}$-module and that we have $\End_{l}(\tilde{V}) \cong \tilde{S}$ as $(\tilde{S}, \tilde{S})$-bimodules. Thus, the $S$-module $V$ descends to $l$ via $\tilde{V}$.

            \item Since $V$ is an endopermutation $kP$-module with $P$ abelian, by \cite[Theorem 7.2]{R96} $V$ has an endosplit $p$-permutation resolution $X_V$. The uniqueness of the lifting $\sigma\colon P \to S^\times$ of the action of $P$ on $S$ with determinant 1 implies the isomorphism class of $V$ is $E$-stable, and thus $X_V$ is an $E$-stable endosplit $p$-permutation resolution. Then $V^* \otimes_k X_V$ is an endosplit (non-$p$-permutation) resolution of $V^* \otimes_k V$, and it follows from \cite[Lemma 1.5(ii)]{HL00} that $\Ind^L_P(V^* \otimes_k X_V)$ is an endosplit resolution of $\Ind^L_P(V^* \otimes_k V),$ where $L := G/O_{p'}(G)$. Since $V$ is capped, $V^*\otimes_k V$ has the trivial $kP$-module as a direct summand, thus $\Ind^L_P(V^* \otimes_k V)$ has the trivial $kL$-module as a direct summand, since $E$ is a $p'$-group. Therefore, $\Ind^L_P(V^* \otimes_k X_V)$ has a direct summand $Y$ which is an endosplit resolution of the trivial $kL$-module. Setting \[X = \Ind^{\hat{L} \times \check{L}}_{\Delta L}(Y),\] we have by Proposition \ref{prop:indcommute} that \[X^* \otimes_{k_\ast \hat{L}} X \cong \Ind^{\hat{L} \times \check{L}}_{\Delta L}(Y^* \otimes_k Y),\] and this is a split complex with homology concentrated in degree 0 isomorphic to $\Ind^{\hat{L} \times \check{L}}_{\Delta L} (k) \cong k_\ast \hat{L}$. Thus $X$ is a Rickard tilting complex of $(k_\ast \hat{L}, k_\ast \hat{L})$-bimodules.

            Since $V$ descends to $\F_p$, $X_V$ descends to $l$ by Theorem \ref{thm:descentofepprs} (note that it remains $E$-stable since Theorem \ref{thm:galoisstablelift} states that descent is unique). By construction, $Y$ also descends to $l$. $k_\ast \hat{L}$ also descends to $l$ since $k\otimes_l l_\ast\hat{L} \cong k_\ast \hat{L}$, so the functor $\Ind^{\hat{L} \times \check{L}}_{\Delta L}$ is well defined, hence $X$ also descends to a chain complex $\tilde{X}$ of $(l_\ast\hat{L}, l_\ast\hat{L})$-bimodules satisfying $\tilde{X}^* \otimes_{l_\ast \hat{L}} \tilde{X} \simeq l_\ast \hat{L}$.
        \end{itemize}

        It follows that $V\otimes_k X$ is a Rickard tilting complex of $(S \otimes_k k_\ast \hat{L}, k_\ast \hat{L})$-bimodules with homology concentrated in degree 0 isomorphic to $V \otimes_k k_\ast \hat{L}$. Then $V \otimes_k X$ is verified to be splendid on \cite[Page 95]{HL00}, completing the proof for the case $G = O_{p',p,p'}(G)$ with abelian Sylow $p$-subgroup. We have shown above that $V \otimes_k X$ descends to a Rickard tilting complex $\tilde{V} \otimes_l \tilde{X}$ of $(\tilde{S} \otimes_l l_\ast \hat{L}, l_\ast \hat{L})$-bimodules with homology concentrated in degree 0 isomorphic to $\tilde{V} \otimes_l l_\ast \hat{L}$. The same proof also shows that $\tilde{X} \otimes_l \tilde{V}$ is splendid. From this we obtain the following.
    \end{construction}

    \begin{theorem}\label{thm:mainthm}
        Let $G$ be a $p$-solvable group satisfying $G = O_{p',p,p'}(G)$ with a abelian Sylow $p$-subgroup $P$. Let $b'$ be a $G$-stable block of $\overline{\F}_pO_{p'}(G)$, and let $b$ be a block idempotent of $\F_p G$ for which $bb' \neq 0$. Let $c$ be the Brauer correspondent of $b$ in $\F_pN_G(P)$. Then $b$ and $c$ are splendidly Rickard equivalent.
    \end{theorem}
    \begin{proof}
        Proposition \ref{prop:structureforreducedcase} implies $b'$ is a block of $\overline{\F}_pG$ with defect group $P$. Let $c'$ be the Brauer correspondent of $b'$, in which case we have $cc' \neq 0$ as well. Then \cite{HL00} proves $b'$ and $c'$ are splendidly Rickard equivalent, and Remark \ref{rem:reduction} demonstrates that the equivalence may be realized over some finite field $k$. Construction \ref{cons:descentofspl} verifies that the construction of \cite{HL00} descends to $k[b] = k[c]$, and \cite[Theorem 6.5(b)]{KL18} implies that $b$ and $c$ are splendidly Rickard equivalent over $\F_p$ (after using Rickard's lifting theorem \cite[Theorem 5.2]{R96} to lift the equivalence of $b$ and $c$ to an absolutely unramified complete discrete valuation ring).
    \end{proof}

    If $G= O_{p',p,p'}(G)$ with abelian Sylow $p$-subgroup, we say a block idempotent $b$ of $\F_pG$ \textit{satisfies condition $(\ast)$} if it arises in the situation of Theorem \ref{thm:mainthm}. That is, there exists a $G$-stable block $b'$ of $\overline{\F}_pO_{p'}(G)$ satisfying $bb' \neq 0$. As an immediate corollary, we obtain another proof that various refinements of the Alperin--McKay conjecture proposed by Isaacs--Navarro in \cite{IN02}, Navarro in \cite{N04}, and Turull in \cite{T13} hold in our setting. We refer the reader to \cite[Page 5]{Bo22} for an overview of these refinements of the Alperin--McKay conjecture. Note that Turull's reformulation, the strongest of the three refinements, was shown by Turull in \cite{T13} for all $p$-solvable groups using separate methods.

    \begin{corollary}
        Refinements (5), (7), and (8) of the Alperin--McKay conjecture, as given in \cite{Bo22}, hold for blocks satisfying condition $(\ast)$.
    \end{corollary}
    \begin{proof}
        This is an immediate corollary of \cite[Theorem 1.4]{Bo22}, after applying unique lifting of splendid Rickard equivalences from $k$ to $\calO$, since splendid Rickard equivalences induce $p$-permutation equivalences (see \cite[Theorem 1.5]{BX07}).
    \end{proof}

    \begin{remark}
        The key obstruction to verifying Kessar--Linckelmann's conjecture for all $p$-solvable groups by proving the reduction step \cite[Proposition 3.1]{HL00} is that \cite[Lemma 3.3]{HL00} may not hold for fields not large enough for $G$. Indeed, a counterexample is as follows. Consider $\F_2S_3 = \F_2(C_3 \rtimes C_2)$, which has two blocks. The non-principal block idempotent $e$ is a block of $\F_2 C_3$. However, in $\F_4 C_3$, $e$ decomposes as $e = e_1 + e_2$, where the two blocks $e_1, e_2$ correspond to the two nontrivial characters of $C_3$, and $C_2$ swaps $e_1$ and $e_2$. It follows that $e$ is a $S_3$-stable block of $\F_2O_{2'}(S_3) = \F_2 C_3$, but $e$ has trivial defect group, not defect group the Sylow 2-subgroup $C_2.$

        This also demonstrates that not all block idempotents of $\F_p G$, where $G$ satisfies $G = O_{p',p,p'}(G)$ with abelian Sylow $p$-subgroup, satisfy condition $(\ast)$. Indeed, with the above notation, there does not exist a $G$-stable block $e'$ of $\overline{\F}_2C_3$ such that $ee' \neq 0$, as $e_1$ and $e_2$ are not $G$-stable.
    \end{remark}

    \textbf{Acknowledgments:} The author thanks Markus Linckelmann for his assistance which helped the author understand certain arguments in \cite{HL00}, and a referee for pointing out a faulty lemma and providing the above counterexample for a previous version of this paper, which had attempted to prove \cite[Proposition 3.1]{HL00} over an arbitrary field of characteristic $p$.

    \bibliography{bib}
    \bibliographystyle{alpha}

\end{document}